\documentclass[12pt]{article}
\usepackage{amsmath,amsfonts,amsthm,amssymb}
\usepackage{epsfig}
\newcommand{\eps}{\varepsilon}

\newtheorem {thm}     {Theorem}

\newtheorem {rmk}  {Remark}

\newtheorem {assump} {Assumption}

\def\({\begin{eqnarray}}
\def\){\end{eqnarray}}
\def\[{\begin{eqnarray*}}
\def\]{\end{eqnarray*}}

\def\R{\mathbb{R}}

\def\d{{\mathrm{d}}}
\def\e{{\mathrm{e}}}
\def\O{{\mathcal{O}}}

\newcommand{\mf}{\bar{f}}

\def\bydef{:=}

\title{Discrete {\em versus} continuous models in evolutionary
dynamics: from simple to
simpler --- and even simpler --- models}

\author{Fabio A. C. C. Chalub\thanks{%
Departamento de Matem\'atica and Centro de Matem\'atica
e Aplica\c c\~oes, Universidade Nova de Lisboa, 
Quinta da Torre, 2829-516, Caparica, Portugal.
e-mail:{\tt chalub@cii.fc.ul.pt}
+351 933 313 096},
Max O. Souza\thanks{%
Departamento de Matem\'atica Aplicada, Universidade Federal
Fluminense, R. M\'ario Santos Braga, s/n, 24020-140, Niter\'oi, RJ, Brasil.
e-mail:{\tt msouza@mat.uff.br}}}

\date{\today}

\begin{document}

\maketitle

\begin{abstract}
There are many different models---both continuous and discrete---used
to describe gene mutation fixation. In particular, the Moran process,
the Kimura equation and the replicator dynamics are all well known
models, that might lead to different conclusions. We present a
discussion of a unified  framework to embrace all these models, in the
large population regime.
\end{abstract}

\section{Introduction}

Real world models need to cover a large range of scales. However,
models that are valid in such a large range are hard to obtain and can
be very complex to analyze. Alternatively, we might use models that focus on
certain scales. Thus, on one hand, we might have a microscopic discrete model
that is derived from first principles while, on the other hand, we
might also have continuous models that are easier to analyze  but are
more phenomenological in nature.

When dealing with many descriptions of the same reality, 
the connection between these various possible descriptions is
an important problem. These connections, in the particular case
of evolutionary game theory for large populations, will be discussed
here after the work in \cite{ChalubSouza}. In particular, we will
present a unified theory that covers the Moran process \cite{Moran},
the mean-field theory description by Kimura \cite{Kimura} and the replicator
dynamics \cite{SigmundHofbauer}.

In order to develop such a theory, we shall proceed as follows:
\begin{enumerate}
\item we prepare a detailed discrete model---the \textsc{the
  microscopic model};
\item we identify suitably scalings and the corresponding {\em
  negligible} (small) parameters; 
\item we formally find a new model where these variables are
set to zero---\textsc{the thermodynamical limit};
\item we prove that the new model is a good approximation
for the previous model, within the given scalings;
\item we study the behaviour of distinguished limiting cases.
\end{enumerate}
It is important to stress that, usually, step three is obtained from
phenomenological framework. Thus, even  formal connections between
discrete and continuous models can be very important in understanding
their relationship. Moreover,  this allows one to solve the simpler
continuous model and thus have the approximate behavior of the
solution of the detailed model. This approach is classical in the
physical sciences where, for instance, continuous mechanics can be
seen as a formal limit of particle dynamics---although the
phenomenological derivation has been obtained much 
earlier~\cite{MarkowichRinghoferSchmeiser}.

In these derivations, the existence of small parameters is generally
natural, but the appropriate scalings are not. 
For example, models of dilute gas given by Boltzmann equation
converge to the Navier-Stokes or Euler equations in fluid
dynamics (depending on the precise scaling given) when 
the re-scaled free mean path is set equal to 
zero~\cite{BardosGolseLevermoreI,BardosGolseLevermoreII,MasmoudiSaint-Raymond,Masmoudi}.  
A similar approach uses kinetic models for modeling cell
movement induced by chemicals (chemotaxis) and when the
cell free mean path is negligible, their solutions is
comparable to the solutions of the Keller-Segel 
model~\cite{CMPS,ChalubKang,HillenOthmer,OthmerHillen}. 
In a different framework, relativistic models for particle
motion have as the non-relativistic limit (i.e., the limit when
typical velocities are small compared to the velocity of
light) the Newtonian physics~\cite{MasmoudiNakanishi,Mauser}, where quantum 
equations converge again
to classical physics when (re-scaled) Planck constant
is very small~\cite{BechoucheMauserSelberg,SparberMarkowich}.

For the Moran Process, it has been recently noticed that the inverse
of the population size is the relevant small parameter;
cf. \cite{TraulsenClaussenHauert,ChalubSouza} for instance. 

The outline of this work is as follows: in section~\ref{Moran}, we
discuss the generalized Moran process. This includes the standard
Moran process as a special case, but addresses also the frequency
dependent case. In section~\ref{Thermo}, we review the scalings and
thermodynamical limits found in \cite{ChalubSouza}. The connection
of some of the thermodynamical limits with the Kimura model is
discussed in section~\ref{Kimura}. After that, we briefly outline some
of the mathematical issues involved in the derivation of the
thermodynamical models. Relationship between these limits and the
Replicator Dynamics is discussed in section~\ref{Replicator}. We then
present a series of numerical simulation to illustrate the theory
discussed and compare results in section~\ref{Tour}. 
Some remarks in more general games, where mixed
strategies are allowed are given in section~\ref{remarks}.

\section{The generalized Moran process}
\label{Moran}

\begin{figure}[h]
\begin{center}
\epsfig{file=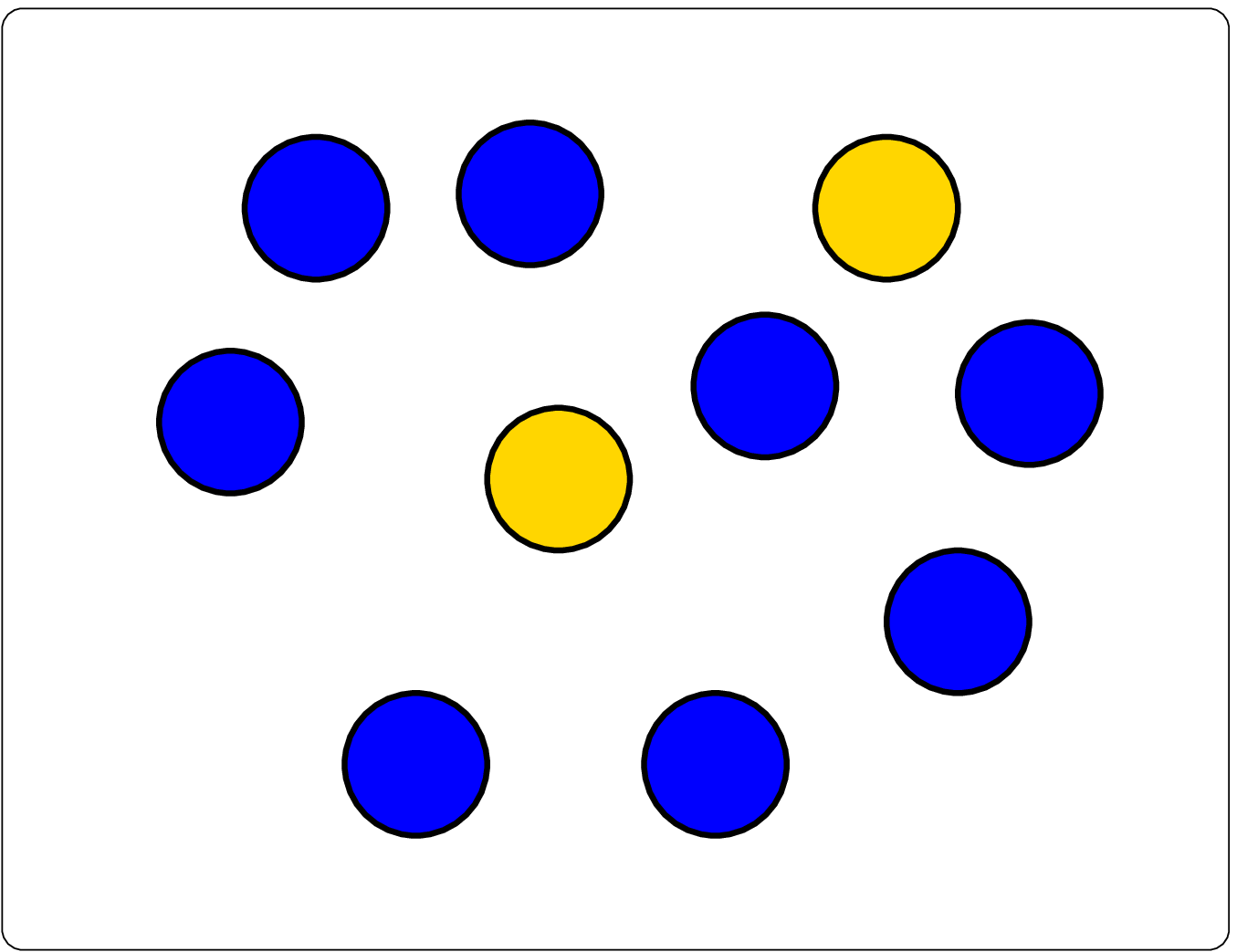,width=0.3\textwidth}
\hspace{.01\textwidth}
\epsfig{file=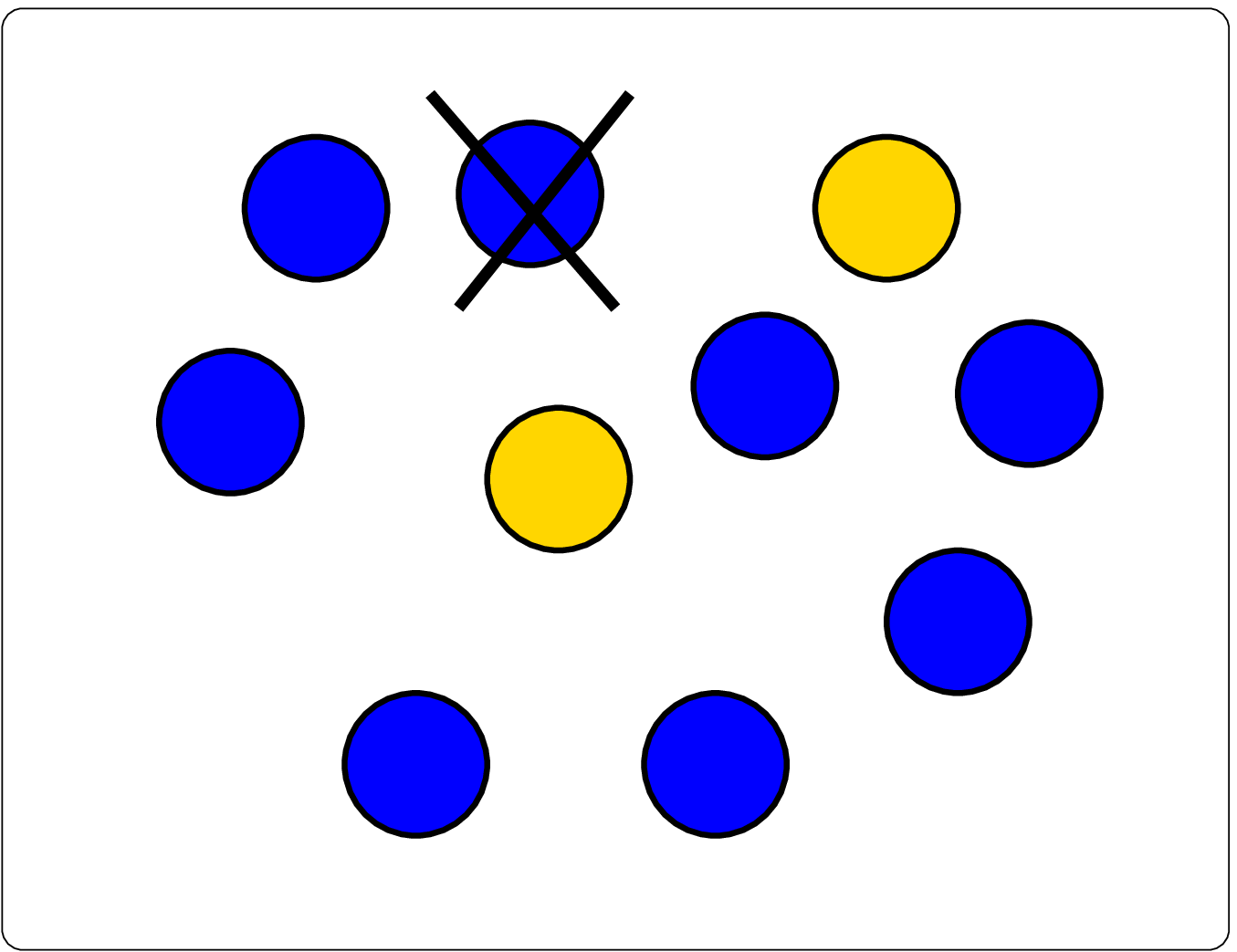,width=0.3\textwidth}
\hspace{.01\textwidth}
\epsfig{file=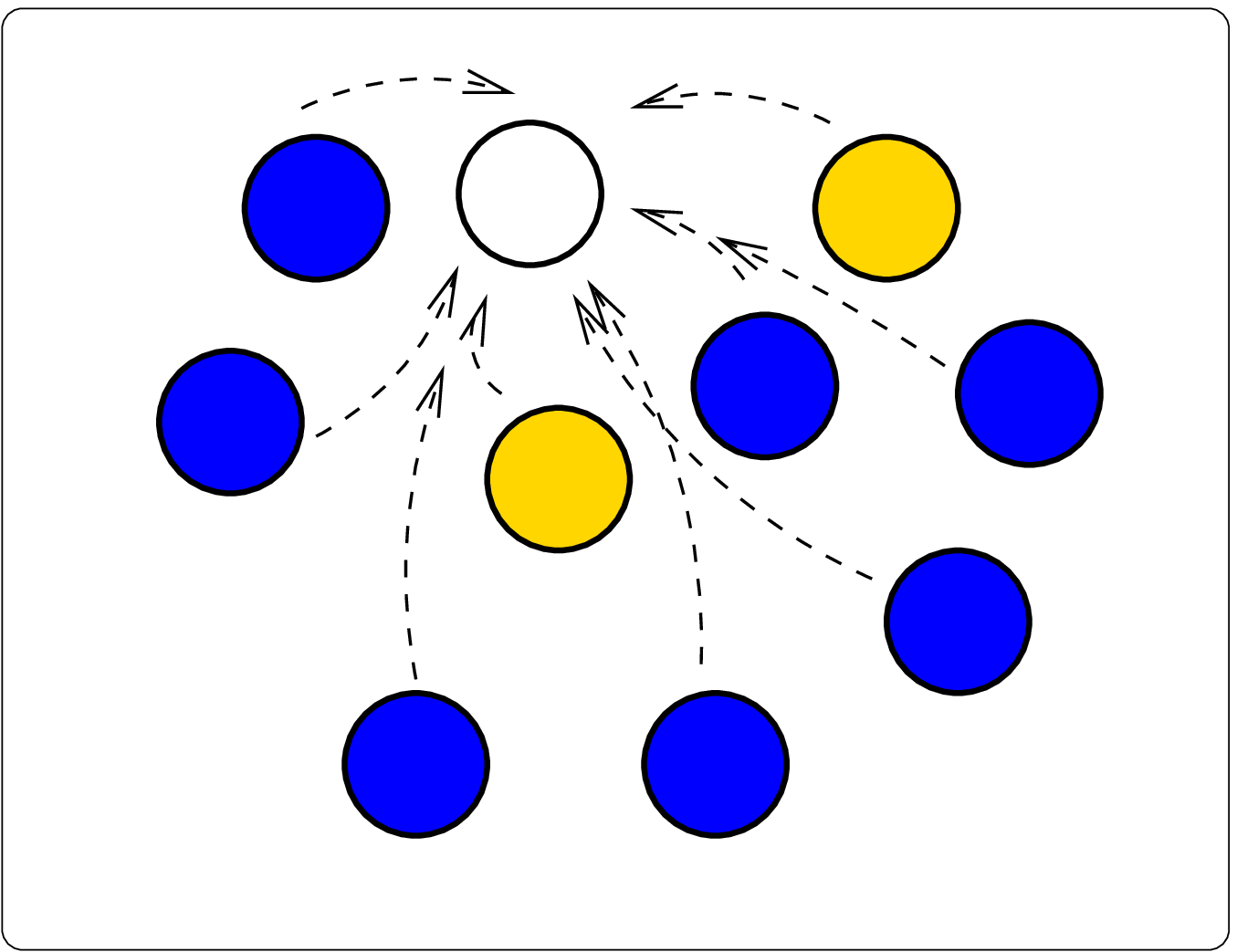,width=0.3\textwidth}
\vspace{.1cm}
\text{(a)\hspace{.29\textwidth} (b)\hspace{.29\textwidth} (c)}
\end{center}
\caption{The Moran process: from a two-types population (a) we chose 
one at random to kill (b) and a second to copy an paste in the place
left by the first, this time proportional to the fitness.
\label{Moran_fig}}
\end{figure}

Consider a population of fixed size $N$, given by
two types of individuals, I and II.
The Moran process is defined by three steps:
\begin{itemize}
\item we choose one of the individuals at random to
be eliminated;
\item all the remaining individuals play a game
given by the pay-off matrix
\begin{center}
\begin{tabular}{c|cc}
&I&II\\
\hline
I&$A$&$B$\\
II&$C$&$D$
\end{tabular}\ ,
\end{center}
and the individual fitness is identified with
the average pay-off. We assume that $C>A>0$ and
$B>D>0$.
This is the only structure of the pay-off matrix that
guarantee the existence of non-trivial stable equilibrium.
This is know in the literature as the ``Hawk-and-Dove'' game.
\item we choose one of the individuals to by
copied with probabilities proportional to the
pay-off.
\end{itemize}
We repeat these steps until a final state is presumably reached.
Intuitively, after a long enough time, all 
individuals will be descendant of a single individual
living at time $t=0$. More precisely, let
$P(t,n,N)$ be the probability that at time $t$ we
have $n$ type I individuals in a total population
of size $N$, and let $c_-(n,N)$ ($c_0(n,N)$, $c_+(n,N)$)
be the transition probability associated to
$n\to n-1$ ($n\to n$ and $n\to n+1$, respectively).
Then, the discrete evolution process is given by
\begin{eqnarray}\label{discrete_evolution}
P(t+\Delta t,n,N)&=&c_+(n-1,N)P(t,n-1,N)\\
\nonumber
&&+c_0(n,N)P(t,n,N)+c_-(n+1,N)P(t,n+1,N)\ .
\end{eqnarray}
Let us introduce the vector
\[
\mathbf{P}(t,N)=(P(t,0,N),P(t,1,N),\ldots,P(t,N,N))^\dagger.
\]
Then the iteration can be written in matrix form as 
\[\mathbf{P}(t+\Delta t,N)=\mathbf{M}\mathbf{P}(t,N),
\]
where $\mathbf{M}$ is a column-stochastic,tridiagonal matrix. It is
also possible to show---cf. \cite{ChalubSouza}---that $1$ is an
eigenvalue of $M$ with multiplicity 
two, with associated eigenvectors given by $\mathbf{e}_1$
and $\mathbf{e}_{N+1}$.
Since the spectrum and its multiplicity is unchanged if $M$ is replaced by
$M^\dagger$ there are two vectors that are kept invariant by
$M^\dagger$. One of them is easily seen to be the vector
$\mathbf{1}=(1,1,\ldots,1,1)^\dagger$. Let $\mathbf{F}$ denote the
remaining one. Then, \cite{ChalubSouza} showed that $\mathbf{F}$ yields
the stationary fixation probability and also that the quantities
\[
\eta_1=\langle \mathbf{1},\mathbf{P}(t,N)\rangle\quad\text{and}\quad
\eta_2=\langle \mathbf{F},\mathbf{P}(t,N)\rangle
\]
are invariants of the Moran process. The former is well known and
it corresponds to the conservation of probability. The latter, however,
seems to be new and it can be interpreted as stating that the correlation
coefficient between a possible state of the Moran process and the
stationary fixation probability must always be the same of the initial
condition. 

We can now prove that
\begin{equation*}
\lim_{\kappa\to\infty}\mathbf{M}^\kappa=
\left(
\begin{matrix}
1&1-F_1&1-F_2&\cdots&1-F_{N-1}& 0\\
0&0&0&\cdots&0&0\\
\vdots&\vdots&\vdots&\ddots&\vdots&\vdots\\
0&0&0&\cdots&0&0\\
0&F_{1}&F_{2}&\cdots&F_{N-1}&1
\end{matrix}
\right)
\end{equation*}
As a direct consequence, for any normalized initial condition
\[
\mathbf{P}_0=(P(0,0,N),P(0,1,N),\cdots,P(0,N,N)),
\]
the final state is 
\[
\lim_{\kappa\to\infty}\mathbf{P}(\kappa\Delta t,N)=
\lim_{\kappa\to\infty}\mathbf{M}^\kappa\mathbf{P}_0=(1-A,0,0,\cdots,A)\ ,
\]
where
\[
A=\sum_{n=0}^NF_nP(0,n,N)
\]
is the {\em fixation probability} associated to the initial condition
$\mathbf{P}_0$. Note that, with probability 1, one of the types will
be fixed. This means, that in the long range, every 
mutation will be either fixed or lost.

\section{Scaling and thermodynamical limits}
\label{Thermo}

The central idea of this section is to find a continuous model that
works as a good approximation of the Moran process, when the
total population is large. This means that we want to find a
continuous model 
for the fraction of mutants in the limit $N\to\infty$. The core
of this process is to define a correct scaling for the time-step
and for the pay-offs. We will show, however, that different
scalings will give different thermodynamical
limits. But only one of these scalings will be able to
capture one essential feature of the discrete process discussed
in the previous section: that genes are always fixed
or lost. In the continuous setting, this means that, as time increases,
the probability distribution should move (diffuse) to
the boundaries.

Let us suppose that (formally) there exists a probability density 
\begin{equation*}
p(t,x)=\lim_{N\to\infty}\frac{P(t,xN,N)}{1/N}=
\lim_{N\to\infty}NP(t,xN,N)\ ,
\end{equation*}
where $x=n/N$.

Let us also suppose that this function $p:\R_+\times[0,1]\to \R$
is sufficiently smooth that we can expand the evolution 
equation~(\ref{discrete_evolution}) as to obtain
\begin{eqnarray}
\label{EqForp}
&&\frac{p(t+\Delta t,x)-p(t,x)}{\Delta t}=\frac{1}{N\Delta t}\left
[\left(c_+^{(1)}+c_0^{(1)}+c_-^{(1)}\right)p
-\left(c_+^{(0)}-c_-^{(0)}\right)\partial_xp\right]\\
\nonumber
&&\qquad
+\frac{1}{N^2\Delta t}\left[\frac{1}{2}\left(c_+^{(2)}+c_0^{(2)}+c_-^{(2)}\right)p-
\left(c_+^{(1)}-c_-^{(1)}\right)\partial_xp
+\frac{1}{2}\left(c_+^{(0)}+c_-^{(0)}\right)\partial_x^2p\right]\\
\nonumber
&&\qquad
+\O\left(\frac{1}{N^3}\right)\ ,
\end{eqnarray}
where $c_*^{(i)}=c_*^{(i)}(x)$, $*=+,0,-$, $i=0,1,2$, are defined by
\begin{eqnarray*}
\label{exp_c+}
c_+\left(xN-1,N\right)&=&c_+(n-1,N)
=c_+^{(0)}+\frac{1}{N}c_+^{(1)}+\frac{1}{2N^2}c_+^{(2)}+\O\left(\frac{1}{N^3}\right)\ ,\\
\label{exp_c0}
c_0(xN,N)&=&c_0(n,N)=c_0^{(0)}+\frac{1}{N}c_0^{(1)}+
\frac{1}{2N^2}c_0^{(2)}+\O\left(\frac{1}{N^3}\right)\ ,\\
\label{exp_c-}
c_-\left(xN+1,N\right)&=&c_-(n-1,N)=
c_-^{(0)}+\frac{1}{N}c_-^{(1)}+\frac{1}{2N^2}c_-^{(2)}+\O\left(\frac{1}{N^3}\right)\ .
\end{eqnarray*}
So far, we made no assumptions on the behavior of the payoffs
$A$, $B$, $C$, and $D$ as $N\to\infty$.  Since the appropriate large
$N$ limit will be attained through a rescaling in time, the fitness
should also rescaled accordingly to preserve the expected amount of
offsprings. In particular, we must have that the fitness approaches
one as $N$ grows large; see discussion in \cite{ChalubSouza}. 
Equivalently, payoffs must also approach one. The order
with they approach unity, however, is a free parameter at this point.

Thus, we impose that
\begin{eqnarray}
\label{Ntoinfty1}
&&\lim_{N\to\infty}(A,B,C,D)=(1,1,1,1)\ ,\\
\label{Ntoinfty2}
&&\lim_{N\to\infty}N^\nu(A-1,B-1,C-1,D-1)=(a,b,c,d)\ ,\quad\nu> 0\ ,
\end{eqnarray}
we find, after a long computation, that
\begin{align*}
\lim_{N\to\infty}N^\nu\left(c_+^{(1)}+c_0^{(1)}+c_-^{(1)}\right)=&
-3x^2(a-b-c+d)-\\&-2x(a-c-2(b-d))+(d-b)\ ,\\
\lim_{N\to\infty}N^\nu\left(c_+^{(0)}-c_-^{(0)}\right)=&
x(1-x)(x(a-c)+(1-x)(b-d))\ .
\end{align*}
The only non-trivial balances, as can be readily observed, are
given by time-steps of order $\Delta t=N^{-(\nu+1)}$, for
$\nu\in(0,1]$. In this case, we have  that
\begin{equation}\label{nodiffusion}
\partial_tp=-\partial_x\left(x(1-x)(x\alpha+(1-x)\beta)p\right)\
,\quad\nu\in(0,1) 
\end{equation}
and
\begin{equation}\label{complete}
\partial_tp=\partial_x^2\left(x(1-x)p\right)-
\partial_x\left(x(1-x)(x\alpha+(1-x)\beta)p\right)\
,\quad\text{if}\quad\nu=1\ , 
\end{equation}
where $\alpha=a-c<0$ and $\beta=b-d>0$. In the particular case where
$\alpha=\beta$ (i.e., when the fitness is independent of the 
particular composition of the population),
the last equation can be shown to be equivalent
  to a celebrated equation in population genetics known as the Kimura
equation~\cite{Kimura}. 

The equations above are supplemented by the following boundary
conditions
\[
\frac{\d}{\d t}\int_0^1p(t,x)\d x=0,\quad\text{for equations
  (\ref{nodiffusion}) and (\ref{complete});}
\]
\[
\frac{\d}{\d t}\int_0^1\psi(x)p(t,x)\d x =0\quad\text{for equation
  (\ref{complete}).}
\]
In the latter condition, $\psi$ satisfies
\begin{equation}
\label{psieq}
\psi''+(\beta+(\alpha-\beta)x)\psi'=0,\quad \psi(0)=0\quad\text{and}
\quad \psi(1)=1.
\end{equation}
The function $\psi(x)$ is the continuous limit of the vector
$\mathbf{F}$ defined in section~$\ref{Moran}$.
\begin{rmk}
For equation (\ref{nodiffusion}), it can be shown that the former condition is
automatically satisfied; hence, we can treat it as a problem with no boundary
conditions. As for equation (\ref{complete}), the degeneracy at the
endpoints with such integral boundary conditions turn it in a very
nonstandard problem in parabolic pdes. We discuss some of the issues
raised by this problem in section~\ref{Math}.
\end{rmk}

\section{The Kimura connection}
\label{Kimura}

Using mean-field Gaussian approximations for the frequency independent
case, Kimura \cite{Kimura} has derived a PDE for the evolution of the
transient fixation probability. which will presumably evolve to a
stationary solution that will then be the standard fixation
probability. This equation is now known as the Kimura equation and
read as follows:
\begin{equation}
\label{Kimura:eq}
\partial_t f = x(1-x)\partial^2_xf +
\gamma x(1-x)\partial_xf,\quad
f(t,0)=0\quad\text{and}\quad
f(t,1)=1.
\end{equation}
The stationary state can be readily found as
\[
f_s(x)=\frac{1-\e^{-\gamma x}}{1-\e^{-\gamma}},
\]
which corresponds to $\psi$ given by (\ref{psieq}), with $\alpha=0$
and $\beta=\gamma$. 

Let $f=\mf(x)+f_s(x)$, then $\mf$ satisfies~(\ref{Kimura:eq}) with
homogeneous boundary conditions. 
In section~\ref{Math}, it will be shown that $p$ can be written as a
sum of a smooth part $q$ with a distributional part with support at the
end points. It will be also shown that $q$ satisfies (\ref{complete})
without boundary conditions.

Now let us assume that $f(t,x)$ is sufficiently smooth in $x$. Then a
straightforward computation shows that
\begin{align*}
\int_0^1\left[x(1-x)\partial^2_x\mf(t,x +
\gamma x(1-x)\partial_x\mf(t,x)\right]q(t,x)\d x&=\\
\int_0^1\mf(t,x)\left[\partial_x^2\left(x(1-x)q(t,x)\right)-
\gamma\partial_x\left(x(1-x)q(t,x)\right) \right]\d x.&
\end{align*}
Thus, (\ref{Kimura:eq}) and (\ref{complete}) with no boundary
conditions are formally adjoints with the appropriate inner product. 

A further relationship between $\mf$ and $q$ should be pointed out,
namely  that, up to a normalising constant, we have
\[
\mf(t,x)=x(1-x)q(t,1-x).
\]
The adjointness discussed above also hold when $\alpha$ is
nonzero. In this case, we have the \textit{generalized Kimura
  equation} given by
\begin{align*}
\partial_t f = x(1-x)\partial^2_xf +
\gamma x(1-x)\left(\beta+(\alpha-\beta)x\right))\partial_xf,&\\
f(t,0)=0\quad\text{and}\quad f(t,1)=1.&
\end{align*}

\section{Mathematical issues}
\label{Math}

There are a number of important questions related to equations
(\ref{complete}) and (\ref{nodiffusion}) given the degeneracy at the
endpoints and the non-standard boundary conditions.

Note that the last two equations will, generally, have very different
qualitative behavior as $t\to\infty$. In particular, we prove the following,
concerning equation~(\ref{complete}):

\begin{thm}\label{main}
\begin{enumerate}
\item For a given $p^0\in L^1([0,1])$, there exists a unique solution
  $p=p(t,x)$  to Equation~(\ref{complete}) of class 
$C^\infty\left(\mathbb{R}^+\times(0,1)\right)$ that satisfies $p(0,x)=p^0(x)$.
\item The solution can be written as
\[
p(t,x)=q(t,x)+a(t)\delta_0+b(t)\delta_1,
\]
where $q\in C^{\infty}(\mathbb{R}^{+}\times [0,1])$ satisfies
(\ref{complete}) without boundary conditions, and we also have
\[
a(t)=\int_0^tq(s,0)\d s\quad\text{and}\quad
b(t)=\int_0^tq(s,1)\d s.
\]
In particular, we have that $p\in C^{\infty}(\mathbb{R}^{+}\times (0,1))$.
\item We also have that
\[
\lim_{t\to\infty}q(t,x)=0\text{ (uniformly)},\quad \lim_{t\to\infty}
a(t)=\pi_0[p^0]\quad \text{and}\quad
 \lim_{t\to\infty} b(t)=\pi_1[p^0],
\]
where  $\pi_0[p^0]=1-\pi_1[p^0]$ and 
{\em the fixation probability} associated to the initial condition
$p^0$ is
\[
\pi_1[p^0]=\frac{\int_0^1\left[\int_y^1p^0(x)\d x\right]
\exp\left(-y^2\frac{\alpha-\beta}{2}-
y\beta\right)\d y}
{\int_0^1\exp\left(-y^2\frac{\alpha-\beta}{2}-
y\beta\right)\d y}\ .
\] 
Note that this means that the solution solution will 'die out' in
the interior and only the Dirac masses in the end points will survive.
\item Write $p^0=a^0\delta_0+b^0\delta_1+q^0 \in L^1([0,1])$ and let
  $\lambda_0$  be the smallest eigenvalue of the associated
  Sturm-Liouville problem (cf. \cite{ChalubSouza}). If, we assum that
$ q^0\in L^2([0,1],x(1-x)\d x)$ and if $||.||_2$ denotes the
  corresponding norm, then we have that
\[
||q(t,.)||_2\leq e^{-\lambda_0t}||q^0(.)||_2.
\]
Moreover, we always have the following $L^1$ bounds:
\begin{enumerate}
\item \[||q(t,.)||_1\leq\e^{-\lambda_0t}||q^0(.)||_1;\]
\item \[\pi_0[p^0]-\e^{-\lambda_0t}||q^0(.)||_1\leq a(t)\leq\pi_0[p^0];\]
\item \[\pi_1[p^0]-\e^{-\lambda_0t}||q^0(.)||_1\leq b(t)\leq\pi_1[p^0].\]
\end{enumerate}
\end{enumerate}
\label{pde_prop}
\end{thm}
It is important to note that equation~(\ref{nodiffusion}) is not a good
long-term approximation for the discrete process in the case of a Hawk-Dove
game, as will see that it presents no diffusion to the boundaries. In
this case, the final state of any  non-trivial initial condition will
be fully determined by the unique  non-trivial equilibrium of the
game, as the following result shows:

\begin{thm}
Consider $p(t,x)\in (L^1\cap H^{-1})([0,1])$ solution of  
Equation~(\ref{nodiffusion}). Then
\[
\lim_{t\to\infty}p(t,x)=\delta_{x^*}\ ,
\]
where $x^*=\beta/(\beta-\alpha)$.
\end{thm}

\begin{proof}
Consider 
\[
\phi(x)=\frac{(x(\alpha-\beta)+\beta)^{\frac{\alpha-\beta}{\alpha\beta}}}
{x^{\frac{1}{\beta}}(1-x)^{-\frac{1}{\alpha}}}\ ,\qquad-\frac{1}{\alpha},\frac{1}{\beta},
\frac{\alpha-\beta}{\alpha\beta}>0\ .
\]
Then, $x(1-x)(x(\alpha-\beta)+\beta)\phi'(x)=-\phi(x)$, which implies
\[
\partial_t\int_0^1p(t,x)\phi(x)\d x=-\int_0^1p(t,x)\phi(x)\d x\ ,
\]
and we conclude that the final state is supported at
$x^*$, the only zero of $\phi(x)$. Using the conservation of mass,
we prove the theorem. 
\end{proof}
\begin{rmk}
\begin{enumerate}
\item For the case of non Hawk-Dove game, i.e., a game only with
  trivial stable equilibrium, then we have 
\[
\lim_{t\to\infty}p(t,x)=c\delta_0+(1-c)\delta_1,
\]
The constant $c$ is directly related to the fixation probability, in
the following sense. Let $\pi_0[p_0]_\epsilon$ be the fixation
probability found with $\alpha$ and $\beta$ replaced by
$\eps^{-1}\alpha$ and $\eps^{-1}\beta$ respectively. 
Then,~\cite{ChalubSouza2006b} show that
\[
c=\pi_0[p^0]+\O(\eps).
\]
Thus, if $\alpha$ and $\beta$ in the original problem are interpreted
as scaled down selection parameters, then equation (\ref{nodiffusion})
yields the same asymptotic behaviour.
\item For initial conditions in $L^1([0,1])$, an
  adaptation of the boundary coupled weak solution developed in
  \cite{ChalubSouza} may be used to show similar results for games
  with or without a non-trivial equilibrium.
\end{enumerate}
\end{rmk}

Equation~(\ref{complete}) is a good approximation for the
discrete case, as can be seen in the following:
\begin{thm}
Let $p_{N,\Delta t}(x,t)$ be the solution of the finite population
dynamics (of population $N$, time step $\Delta t=1/N^2$), with initial
conditions given by $p^0_N(x)=p^0(x)$, $x=0,1/N,2/N,\cdots,1$, for
$p^0 \in L^1_+([0,1])$. Assume also that $(A-1,B-1,C-1,D-1)=
1/N(a,b,c,d)+\O(1/N^2)$.
Let $p(t,x)$ be the solution of equation~(\ref{complete}),
with initial condition given by $p^0(x)$.
If we write $p_i^n$ for the $i$-th component of $p_{N,\Delta t}(x,t)$
in the $n$-th iteration, we have, for any $t^{*}>0$, that
\[
\lim_{N\to\infty}  p^{tN^2}_{xN}
=p(t,x),\quad x\in[0,1],\quad t\in[0,t^{*}].
\]
\end{thm}

Equation~(\ref{nodiffusion}) is however a good approximation
of~(\ref{complete}) for intermediate times and strong selection.

In fact,
\begin{thm}
Consider $(\alpha',\beta')=\eps(\alpha,\beta)$
and $t'=\eps^{-1}t$. Then, in the limit $\eps\to 0$,
we have that the regular part of the solution $q_\eps$ of the re-scaled
equation~(\ref{complete}) converges to the solution of
equation~(\ref{nodiffusion}) 
in $L^2([0,T]\times[0,1],\d t\otimes x(1-x)\d x)$, if the initial
condition is in $H^1([0,1],x(1-x)\d x)$.
\end{thm}.

\begin{proof}
Dropping $'$, and having in mind Theorem~\ref{main}
we re-write Equation~(\ref{complete})
as
\begin{equation}\label{complete_q}
\partial_tq_\eps=\eps\partial_x^2\left(x(1-x)q_\eps\right)
-\partial_x\left(x(1-x)(x(\alpha-\beta)+\beta)q_\eps\right)
\end{equation}
and then we have the {\em a priori} estimate
\begin{eqnarray*}
&&\frac{1}{2}\int_0^1 x(1-x)q_\eps^2\d x\\
&&\qquad=-\eps\int_0^1\left(\partial_x\left(x(1-x)q_\eps\right)\right)^2\d x
+\frac{1}{2}\int_0^1(x\alpha+(1-x)\beta)\partial_x\left(\left(x(1-x)q_\eps\right)^2\right)\\
&&\qquad\le\frac{\beta-\alpha}{8}\int_0^1x(1-x)q_\eps^2\d x\ .
\end{eqnarray*}

We differentiate equation~(\ref{complete_q}) with respect to $t$,
proceed as above to find the estimate
\[
\int_0^1x(1-x)\left(\partial_t q_\eps\right)^2\d x\le \Phi_1(t)\ ,
\]
for an $\eps$-independent function $\Phi_1$. In order to find an $\eps$-independent
bound for $\int_0^1x(1-x)\left(\partial_xq_\eps\right)\d x$, first we prove
\[
\frac{1}{2}\partial_t\int_0^1q_\eps^2\d x\le
\frac{1}{2}\int_0^1\left[x(1-x)(x(\alpha-\beta)+\beta)-\eps(1-2x)\right]
\partial_x q_\eps^2\d x\le C\int_0^1q_\eps^2\d x\ ,
\]
and this implies an {\em a priori} bound for $\int_0^1q_\eps^2\d x$.
Then, note that
\begin{eqnarray*}
&&\frac{1}{2}\partial_t\int_0^1x(1-x)\left(\partial_xq_\eps\right)^2\d x\le\\
&&\qquad\int_0^1\left[-\frac{3}{2}(\alpha-\beta)x(1-x)-(1-2x)(x\alpha+(1-x)\beta)\right]
\left(\partial_xq_\eps\right)^2\d x\\
&&\qquad\quad+\int_0^1\left[(x\alpha+(1-x)\beta)-(\alpha-\beta)(1-2x)\right]
x(1-x)\partial_xq_\eps^2\d x\\
&&\qquad\le C_1\int_0^1x(1-x)\left(\partial_xq_\eps\right)^2\d x
+C_2\int_0^1q_\eps^2\d x\ .
\end{eqnarray*}
We conclude an {\em a priori} bound for 
$\int_0^1x(1-x)\left(\partial_x q_\eps\right)^2\d x$
and then from Rellich's theorem, we know that $\int_0^1x(1-x)\left(q_\eps\right)^2\d x$
is in a compact set of $L^2([0,T]\times[0,1])$. This proves the theorem.
\end{proof}

\begin{rmk}
Equation (\ref{complete}) and (\ref{nodiffusion}) have a very
important difference, even in the case where their asymptotic
behaviour is the same. For equation (\ref{complete}) the
Diracs at the endpoints appear at time $t=0^+$, while for
(\ref{nodiffusion}) this is only attained at $t=\infty$. Thus, we have
the unusual situation that, at the endpoints,  the parabolic problem
is more singular  than the hyperbolic associated problem.
\end{rmk}

\section{The replicator dynamics connection}
\label{Replicator}

The replicator dynamics models the evolution of the fraction
of a given type of individuals in a infinite population
framework. For a pay-off matrix given by
\begin{equation}
\left(
\begin{matrix}
a&b\\
c&d
\end{matrix}
\right)\ ,
\end{equation} 
in its simplest form
the replicator dynamics reads as following
\begin{equation}\label{replicatorODE}
\dot X=X(1-X)(X(\alpha-\beta)+\beta)\ .
\end{equation}
Equation~(\ref{nodiffusion}) can be written as 
\begin{align*}
\partial_tp&+x(1-x)(\beta+(\alpha-\beta)x)\partial_xp+\\
&+\left(\beta + 2(\alpha-2\beta)x - 3(\alpha-\beta)x^2\right)p=0
\end{align*}
Its characteristics are given by
\begin{align*}
\frac{\d t}{\d s} &= 1,\\
\frac{\d x}{\d s} &= x(1-x)(\beta+(\alpha-\beta)x),\\
\frac{\d z}{\d s} &= -\left(\beta + 2(\alpha-2\beta)x -
3(\alpha-\beta)x^2 \right)z.
\end{align*}
The projected characteristics in the $x\times t$ plane are given  by
\[
\frac{\d x}{\d t}=x(1-x)(\beta+(\alpha-\beta)x),
\]
which is just (\ref{replicatorODE}).

For smooth solutions, one can then write the solution to
(\ref{nodiffusion})---as done in \cite{ChalubSouza2006b}---as
\begin{equation}
p(t,x)=a^0\delta_0+b^0\delta_1
+q^0(\Phi_{-t}(x))
\left|\frac{\beta+(\alpha-\beta)\Phi_{-t}(x)}{\beta+(\alpha-\beta)x}\right|
\frac{\Phi_{-t}(x)\left(1-\Phi_{-t}(x)\right)}{x(1-x)},
\label{soln:nodiffusion}
\end{equation}
where $\Phi_t(x)$ is the flow map of (\ref{replicatorODE}).

Notice that, when $\alpha-\beta\not=0$, the first order term does not
represents a pure drift, but also a dampening 
(enhancing) for $\alpha>\beta$
($\alpha<\beta$, respectively).

Thus, equation (\ref{nodiffusion}) can be seen as an Eulerian
representation of a quantity associated to the probability density
evolution, but not to the probability density itself. If we
let $q(t,x)=p(t,x)-a^0\delta_0-b^0\delta_1$, we see, from
(\ref{soln:nodiffusion}), that the Lagrangian transported
  by the replicator flow is 
\[
u(t,x)=x(1-x)(\beta+(\alpha-\beta)x)q(t,x).
\]
Thus, (\ref{replicatorODE}) can be see as a Langragian representation
of $u$, once the initial probability distribution is given. Since, we
can recover $q$ from $u$, and hence can recover $p$, we have that,
when  there is no diffusion, solutions to (\ref{replicatorODE})
together with initial probability distribution are equivalent to
(\ref{nodiffusion}). 

An interesting question is to quantify how good is the dynamics given
by equation (\ref{nodiffusion})---or equation
(\ref{replicatorODE}) for that matter---as an approximation to the
dynamics of (\ref{complete}) in the case of small diffusion, i.e.,
strong selection. Besides the results already alluded to in
section~\ref{Math}, the following results have also been shown in
\cite{ChalubSouza2006b}:
\begin{enumerate}
\item For games without a non-trivial stable equilibrium, we have  that the
  dynamics of $p$ is well approximated by solutions of
  (\ref{nodiffusion}) over a long time  modulated by an envelope on
  a slow timescale.
\item For games with a non-trivial stable equilibrium, the above holds away
  of such an equilibrium. Near the equilibrium, we have a balance of
  diffusive and selective effects. This prevents the Dirac formation
  at the equilibrium point. 
\item Combining the remarks above, we have, for Hawk-Dove games, that
a non trivial
initial distribution (i.e., that is not peaked at the endpoints)
 tends to peak at the interior equilibrium, and that such a peak takes
 a long time to die out. For an example see figure 11 in \cite{ChalubSouza}
\end{enumerate}

\section{A numerical tour}
\label{Tour}

For a compairison of the discrete and continuos models, as well as an
extensive ensemble of simulations for (\ref{complete}), the
reader is referred to \cite{ChalubSouza}.

Here, we shall focus on compairing the solutions to (\ref{complete})
with solutions to (\ref{nodiffusion}). 
We present two sets of simulations of (\ref{complete});  with
large $\beta$ and large $\alpha$. We then compared the solutions to
(\ref{nodiffusion}) with rescaled time and coefficients. We also
plotted the position with the peak---with rescaled height---with the
peak of the solution to (\ref{complete}). For display conveninence, we
have omitted the very ends of the interval and plotted $\Delta x p$
instead of $p$.
\begin{figure}[htbp]
\begin{center}
\parbox{0.45\textwidth}{\epsfig{file=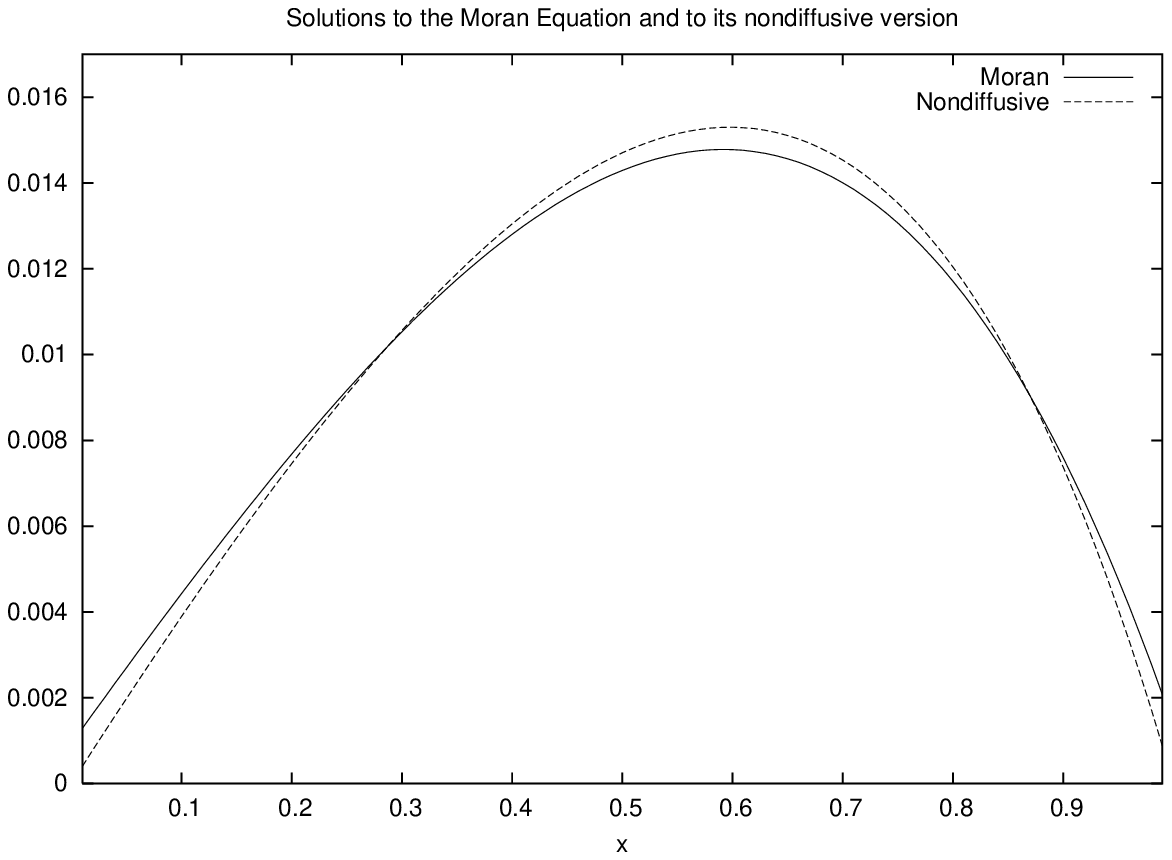,width=0.4\textwidth}}
\quad
\parbox{0.45\textwidth}{\epsfig{file=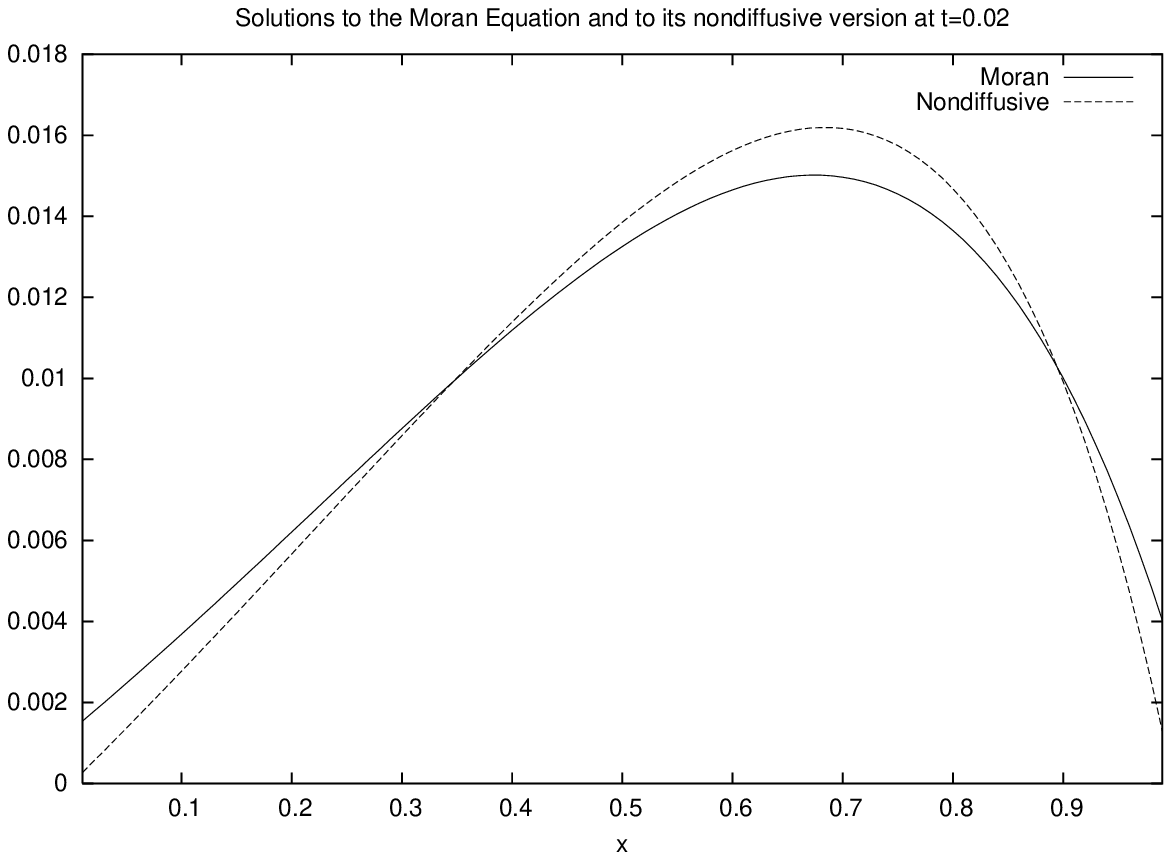,width=0.4\textwidth}}
\caption{Solutions to equation~(\ref{complete}), labled as Moran, and
  to (\ref{nodiffusion}) labeled as Nondiffusive in the frequency
  independent case, with $\alpha=\beta=20$ and initial condition $p^0(x)=x(1-x)/6$.} 
\label{tfig1}
\end{center}
\end{figure}
\begin{figure}[htbp]
\begin{center}
\epsfig{file=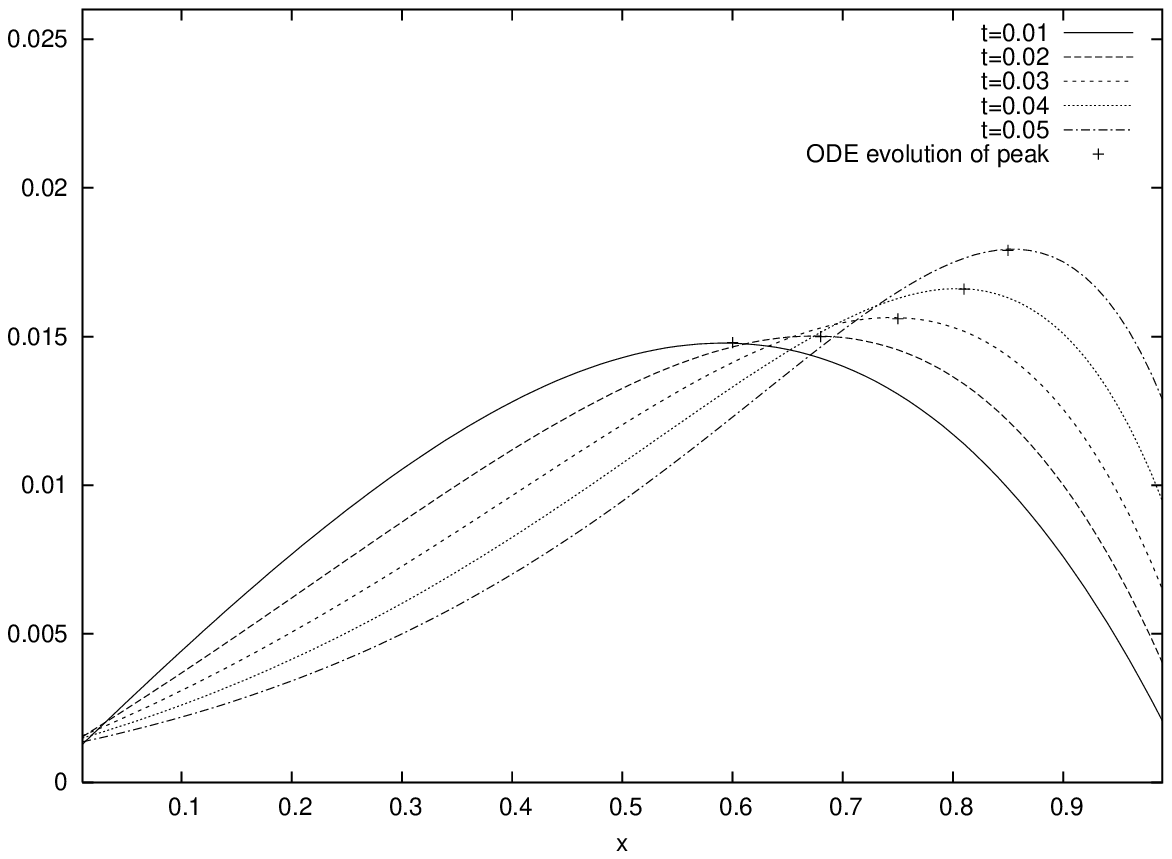,width=0.5\textwidth}
\caption{Evolution of solutions to (\ref{complete}) together with the peaks
  given by solutions to (\ref{nodiffusion}) plotted as points with
  rescaled height for a convenient display. Same paramters and initial
  condition of figure~\ref{tfig1}}
\label{tfig2}
\end{center}
\end{figure}

\bigskip

\begin{figure}[htbp]
\begin{center}
\parbox{0.45\textwidth}{\epsfig{file=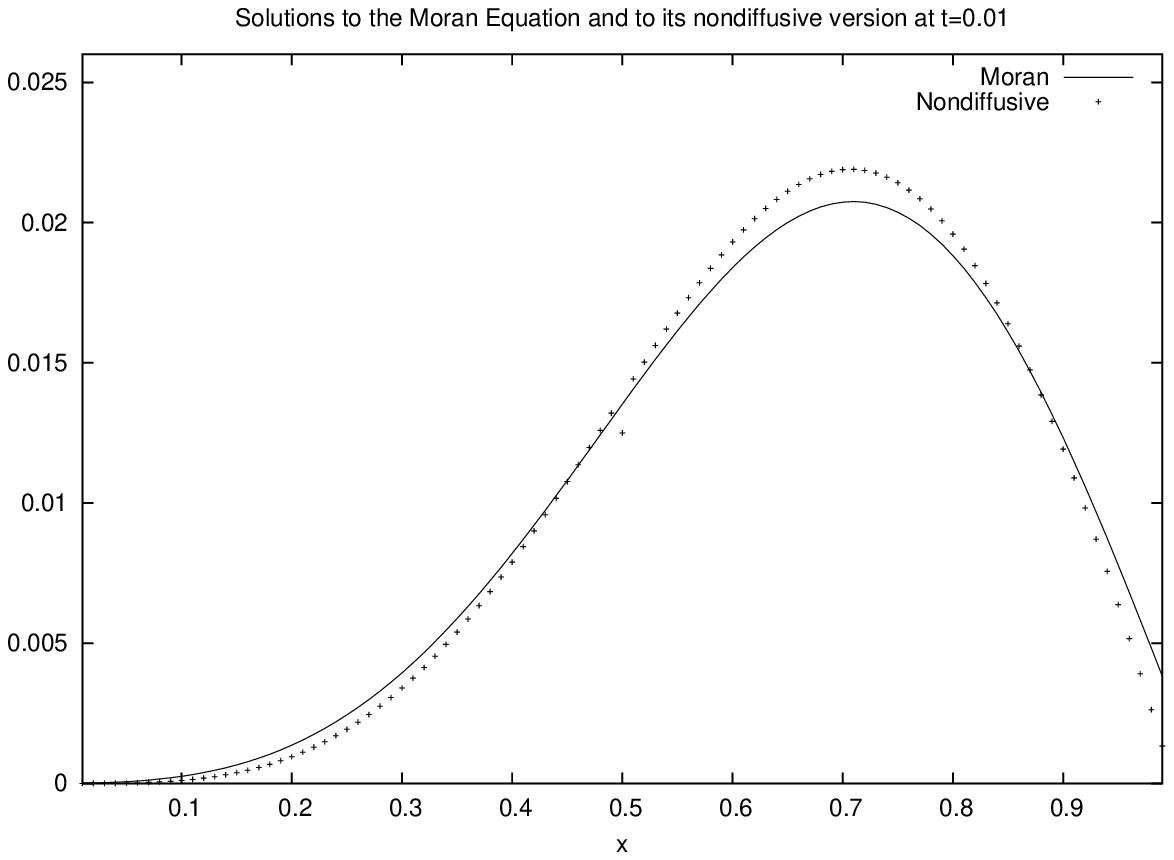,width=0.4\textwidth}}
\quad
\parbox{0.45\textwidth}{\epsfig{file=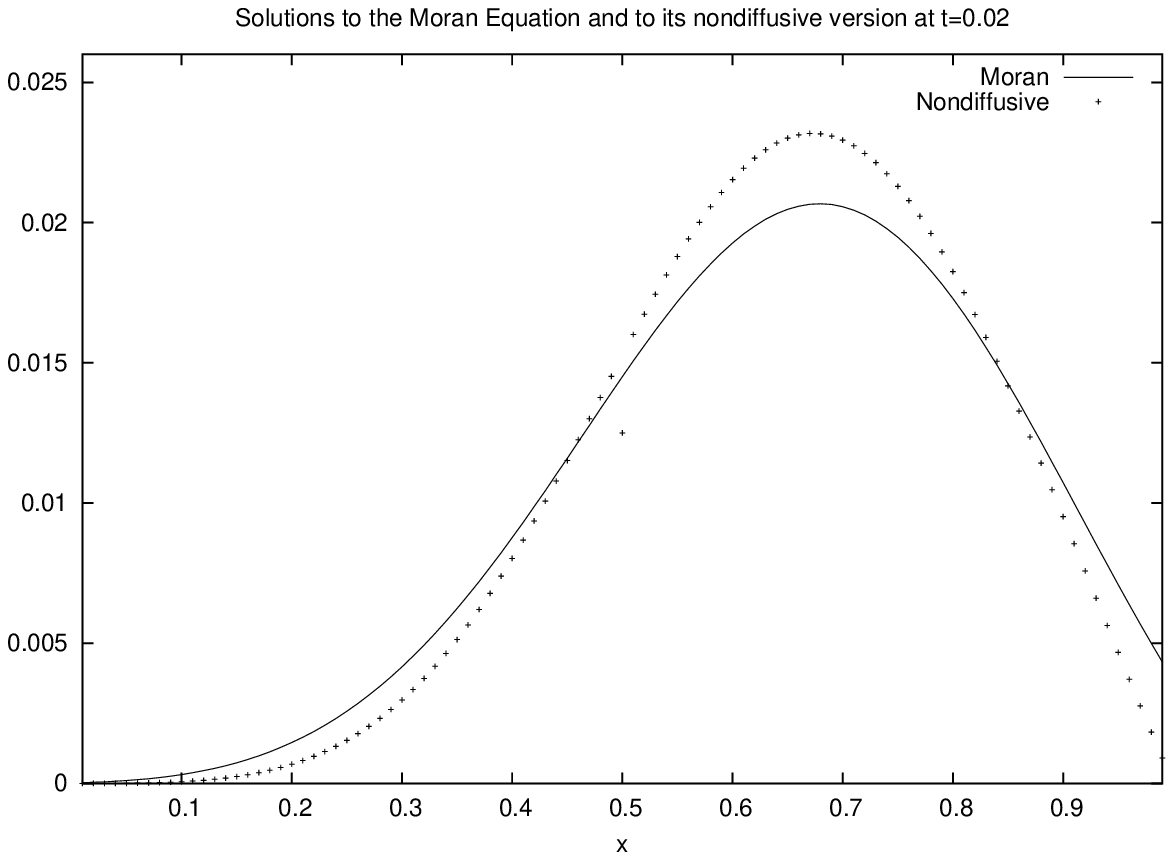,width=0.4\textwidth}}
\caption{Solutions to equation~(\ref{complete}), labled as Moran, and
  to (\ref{nodiffusion}) labeled as Nondiffusive in the frequency
  independent case, with $\alpha=-20$ and $\beta=20$ and initial
  condition $p^0(x)=20x^3(1-x)$.}  
\label{tfig3}
\end{center}
\end{figure}
\begin{figure}[htbp]
\begin{center}
\parbox{0.45\textwidth}{\epsfig{file=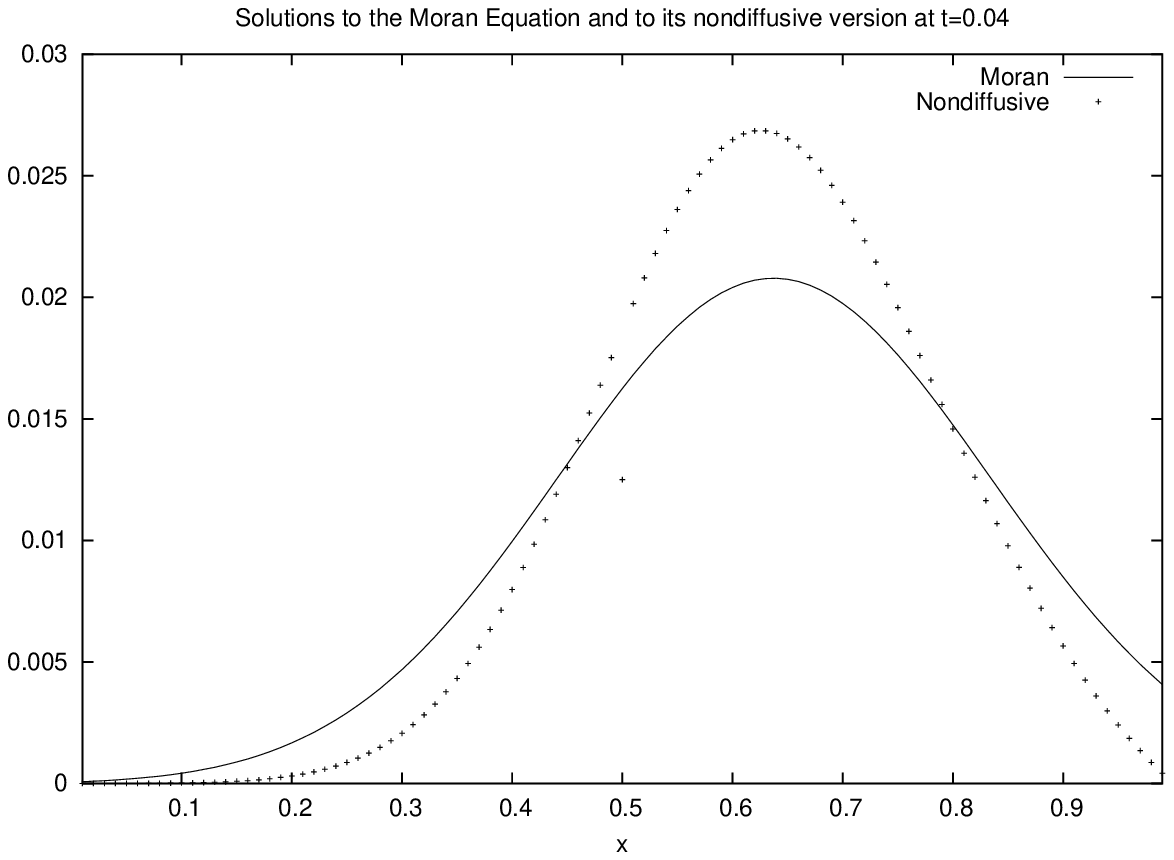,width=0.4\textwidth}}
\quad
\parbox{0.45\textwidth}{\epsfig{file=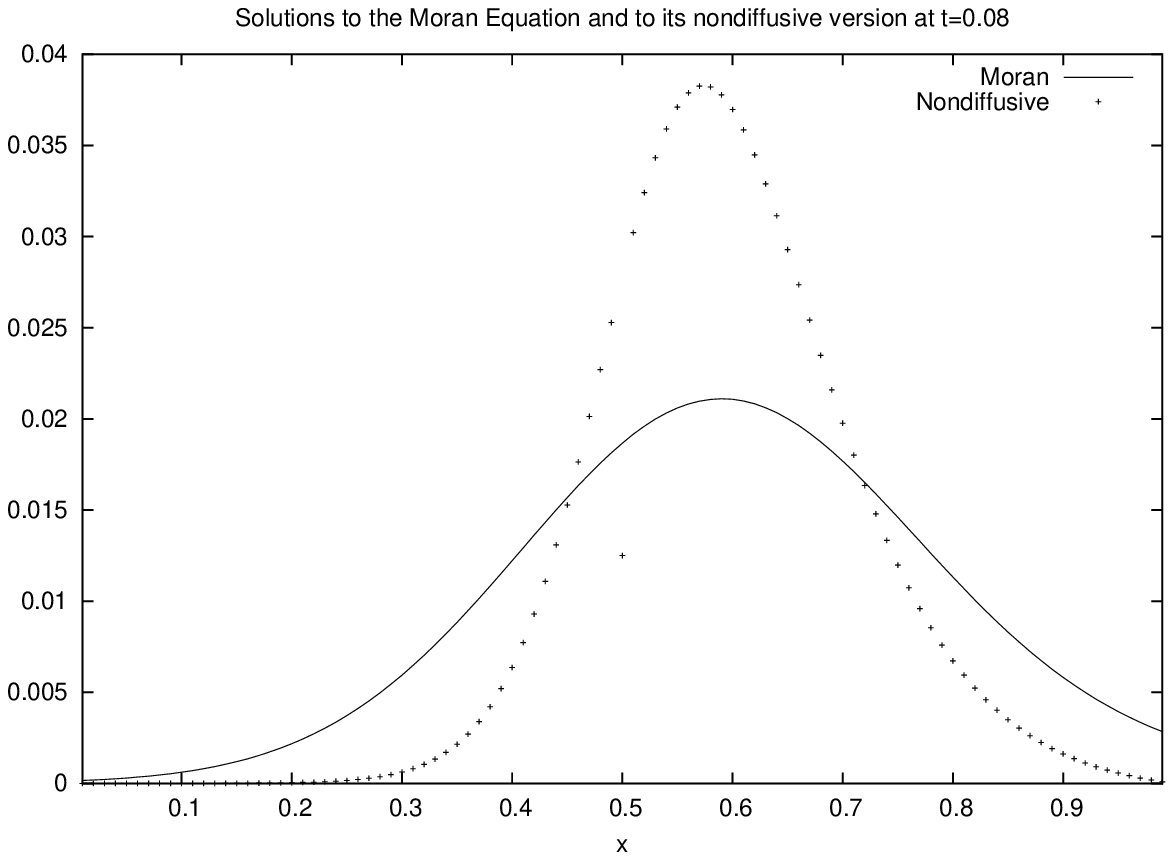,width=0.4\textwidth}}
\caption{Solutions to equation~(\ref{complete}), labled as Moran, and
  to (\ref{nodiffusion}) labeled as Nondiffusive in the frequency
  independent case, with $\alpha=-20$ and $\beta=20$ and initial
  condition $p^0(x)=20x^3(1-x)$.}  
\label{tfig4}
\end{center}
\end{figure}
\begin{figure}[htbp]
\begin{center}
\epsfig{file=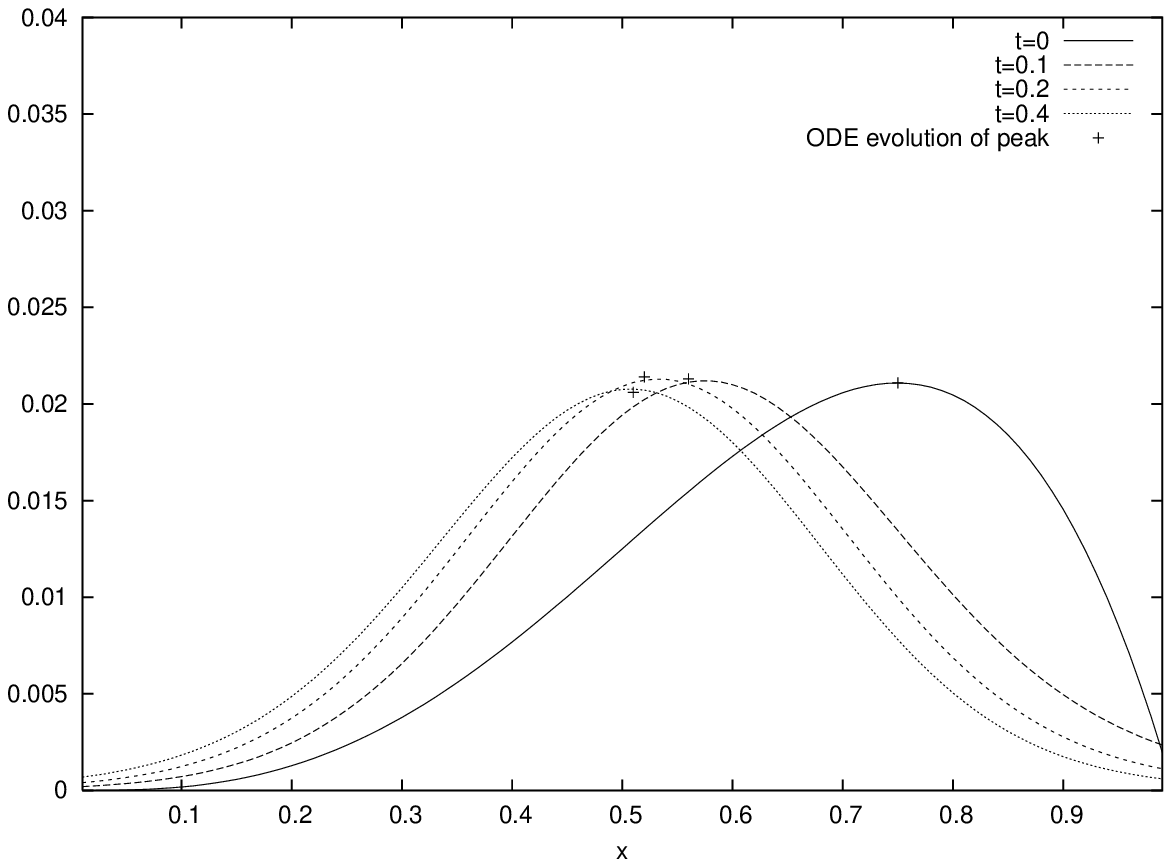,width=0.5\textwidth}
\caption{Evolution of solutions to (\ref{complete}) together with the peaks
  given by solutions to (\ref{nodiffusion}) plotted as points with
  rescaled height for a convenient display. Same paramters and initial
  condition of figures~\ref{tfig3} and~\ref{tfig4}}
\label{tfig5}
\end{center}
\end{figure}

\section{Further remarks}
\label{remarks}

The analogy between the Moran process for finite populations
and the replicator dynamics can be taken further. More
precisely, suppose that the individuals taking part in the
Moran process do not play only pure strategies as in the above
analysis, but are allowed to play mixed strategies. In particular,
let us suppose that the game involves two kind of strategists,
$E_{\theta_1}$ and $E_{\theta_2}$, where an $E_{\theta}$-strategist means that
he/she plays pure strategy I with probability $\theta$ and II with
probability $1-\theta$. Then, the pay-off matrix is given by
\begin{center}
\begin{tabular}{c|cc}
&$E_{\theta_1}$&$E_{\theta_2}$\\
\hline
$E_{\theta_1}$&$\widetilde A$&$\widetilde B$\\
$E_{\theta_2}$&$\widetilde C$&$\widetilde D$
\end{tabular}\ ,
\end{center}
where
\begin{eqnarray*}
\widetilde A&\bydef&\theta_1^2A+\theta_1(1-\theta_1)(B+C)+(1-\theta_1)^2D\ ,\\
\widetilde B&\bydef&\theta_1\theta_2A+\theta_1(1-\theta_2)B+(1-\theta_1)\theta_2C+(1-\theta_1)(1-\theta_2)D\ ,\\
\widetilde C&\bydef&\theta_1\theta_2A+(1-\theta_1)\theta_2B+\theta_1(1-\theta_2)C+(1-\theta_1)(1-\theta_2)D\ ,\\
\widetilde D&\bydef&\theta_2^2A+\theta_2(1-\theta_2)(B+C)+(1-\theta_2)^2D\ .
\end{eqnarray*}
The associated thermodynamical limit is given by
\[
\partial_tp=\partial_x^2\left(x(1-x)p\right)-\partial_x\left(x(1-x)
(x(\theta_1-\theta_2)^2(\alpha-\beta)+(\theta_1-\theta_2)(\theta_2\alpha+(1-\theta_2)\beta p)\right)\ .
\]
The final state is given by $p^\infty=\pi_0[p^0]\delta_0+\pi_1[p^0]\delta_1$,
where $\pi_0[p^0]=1-\pi_1[p^0]$ and
the fixation probability $\pi_1[p^0]$ is given by
\[
\pi_1[p^0]=\frac{\int_0^1\int_0^x p^0(x)F_{(\theta_1,\theta_2)}(y)\d y\d x}
{\int_0^1F_{(\theta_1,\theta_2)}(y)\d y}\ ,
\]
where
\[
F_{(\theta_1,\theta_2)}(y)\bydef\exp\left(-y^2(\theta_1-\theta_2)^2\frac{\alpha-\beta}{2}
-y(\theta_1-\theta_2)(\theta_2\alpha+(1-\theta_2)\beta)\right)\ .
\]
Note that the neutral case (i.e., when the two types of individuals
are of the same kind) is given by $\theta_1=\theta_2$, and in this case
the governing equation is purely diffusive and fixation probability
associated to a given initial state is simply given by
\begin{equation}\label{neutral}
\pi_1^\mathrm{N}[p^0]=\int_0^1xp^0(x)\d x\ .
\end{equation}
We say that an $E_{\theta_2}$ strategist dominates an $E_{\theta_1}$ strategist
($E_{\theta_2}\succ E_{\theta_1}$) if the fixation probability of the first type,
for any non-trivial 
initial condition is smaller that the one in neutral
case given by equation~(\ref{neutral}). With this definition, we can
prove the following theorem:
\begin{thm}
$E_{\theta_2}\succ E_{\theta_1}$ if and only if the flow of the replicator dynamics
is such that $\theta_1\longrightarrow \theta_2$.
\end{thm}

As a simple corollary, we have that if $\theta^*=\beta/(\beta-\alpha)\in(0,1)$
(the ESS of the game), then $E_{\theta^*}\succ E_\theta$, $\forall \theta\not=\theta^*$.
This shows that an individual playing a mixed strategy with probabilities
equal to the one of the game's ESS is better equipped to win any
context. But, as we saw in the previous sections, for populations
of pure strategists we can not expect an stable mixture (even in
fractions equivalent to the game's ESS) to evolve.

\bibliography{ref}
\bibliographystyle{abbrv}

\end{document}